\documentclass[11pt]{article}
\usepackage{amsmath,amssymb,latexsym}
\newtheorem{theo}{Theorem}
\newtheorem{coro}{Corollary}
\newtheorem{ex}{Example}

\newtheorem{prop}{Proposition}

\makeatletter\@addtoreset{equation}{section}\makeatother
\begin{document}
\title{Symmetry and dependence properties within 
a semiparametric family of bivariate copulas}
\author{C\'ecile Amblard$^{1,2}$ \& St\'ephane Girard$^{3}$}
\date{
\small
$^{1}$ LabSAD, Universit\'e Grenoble 2,
BP 47, 38040 Grenoble Cedex 9, France.\\
T\'el : (33) 4 76 82 58 26,  Fax : (33) 4 76 82 56 65,
E-mail : Cecile.Amblard@upmf-grenoble.fr\\
$^{2}$ CRM, Universit\'e de Montr\'eal,
C.P. 6128, succursale centre-ville,\\  Montr\'eal (Qu\'ebec) H3C 3J7, Canada.\\
$^{3}$ INRIA Rh\^one-Alpes, projet {\sc is2}, \\
ZIRST, 655, avenue de l'Europe, Montbonnot, 38334 Saint-Ismier Cedex, France.\\
T\'el : (33) 4 76 61 54 21, Fax : (33) 4 76 61 52 52, E-mail: Stephane.Girard@inrialpes.fr
}
\maketitle

\begin{abstract}
\noindent In this paper, we study a semiparametric family of bivariate copulas.
The family is generated by an univariate function, determining the symmetry
(radial symmetry, joint symmetry) and dependence property (quadrant dependence,
total positivity, ...) of the copulas. We provide bounds on different measures
of association (such as Kendall's Tau, Spearman's Rho) for this family
and several choices of generating functions allowing to reach these bounds.\\



\end{abstract}

                   \section{Introduction}
The theory of copulas provides a relevant tool to build multivariate probability laws, from fixed margins and required degree of dependence. From Sklar's Theorem \cite{SKLAR}, the dependence properties of a continuous multivariate distribution $H$ can be entirely summarized, independently of its margins, by a copula, uniquely associated with $H$. Several families of copulas, such as Archimedian copulas \cite{GENEST} or copulas with polynomial sections \cite{QUESADA,NELSEN97} have been proposed. In \cite{NELSEN97}, the authors point out that the copulas with quadratic section proposed in \cite{QUESADA} are not able to modelize large dependences. 
Then, they introduce copulas with cubic sections and conclude that copulas with higher order polynomial sections would increase the dependence degrees but simultaneously the complexity of the model. 
In \cite{HUE}, polynomial-type extensions of the Farlie-Gumbel-Morgenstern bivariate distributions \cite{FAR} are studied. It is shown that 
the positive correlation between the marginal distributions can be increased up to $\approx 0.39$ but that the maximal negative correlation remains
$-{1/3}$.  
Moreover, it has been remarked that dependence degrees arbitrarily close to
$\pm 1$ cannot be obtained with polynomial functions of fixed degree~\cite{HUE2}.
Thus, we propose to give up the polynomial form to work with a semiparametric family of copulas. The induced parametric families of copulas are 
generated as simply as Archimedian copulas, that is by an univariate function.
Furthermore, dependence properties of copulas with polynomial sections are
preserved and the dependence degree can be increased without significantly
complexifying the model. Note that, in \cite{FER}, a class of symmetric bivariate copulas with wide correlation coefficients range is introduced, but these copulas are less convenient to perform classical calculations on a probability law than the copulas studied here.\\
\noindent In Section \ref{secdef}, the semiparametric family is defined and its basic properties are
derived. Symmetry properties are established in Section \ref{secsym}. Properties of three
measures of association are investigated in Section \ref{secdep}. Section \ref{secdep2} is devoted to
the dependence structure and the dependence ordering of the family.

            \section{Definition and basic properties}
\label{secdef}

\noindent Throughout this paper, we note $I=[0,1]$. A bivariate copula defined on the unit square $I^2$ is a bivariate cumulative distribution function with univariate uniform $I$ margins. Equivalently, $C$ must satisfy the following properties :
\begin{description}
\item [(P1)] $C(u,0)=C(0,v)=0$, $\forall (u,v)\in I^2$,
\item [(P2)] $C(u,1)=u$ and $C(1,v)=v$, $\forall (u,v)\in I^2$,
\item [(P3)] $\Delta(u_1, u_2, v_1, v_2)=C(u_2,v_2)-C(u_2,v_1)-C(u_1,v_2)+C(u_1,v_1)\geq 0$, $\forall (u_1, u_2, v_1, v_2)\in I^4$, such that $u_1\leq u_2$ and $v_1\leq v_2$.
\end{description}
Let us recall that, from Sklar's Theorem, any 
bivariate  distribution with cumulative distribution function $H$
and marginal cumulative distribution functions $F$ and $G$ can be written
$H(x,y)=C(F(x),G(y))$, where $C$ is a copula. This result justifies the use of copulas for building bivariate distributions.\\
\noindent We consider the semiparametric family of functions defined on $I^2$ by :
\begin{equation}
\label{defcop}
C_\theta(u,v)=uv+\theta \phi(u)\phi(v), \;\;\theta\in[-1,1],
\end{equation}
where $\phi$ is a function on $I$. This family was first introduced in \cite{LLALENA}, chapter 3, and is a particular case of Farlie's family introduced in \cite{FAR}. \\
Let us note first that, the independent copula $C_0(u,v)=uv$ belongs to any parametric family $\{C_\theta\}$ generated by a function $\phi$.
Second, the functions $\phi$ and $(-\phi)$ clearly define the same function $C_\theta$.
The following theorem, very similar to Theorem 2.2 in \cite{QUESADA}, dedicated to bivariate copulas with quadratic sections, gives sufficient and necessary conditions on $\phi$ to ensure that $C_\theta$ is a copula. 
\begin{theo}
\label{theocop}
$\phi$ generates a parametric family of copulas $\{C_\theta,\; \theta\in[-1,1]\}$
 if and only if it satisfies the following conditions~:
\begin{description}
\item [(a)] $\phi(0)=\phi(1)=0$,
\item [(b)] $\phi$ satisfies the Lipschitz condition :
$|\phi(x)-\phi(y)|\leq |x-y|, \;\; \forall (x,y)\in I^2$.
\end{description}
Furthermore, $C_\theta$ is absolutely continuous.
\end{theo}
\noindent {\bf Proof: } It is clear that {\bf (P1)}$\Leftrightarrow
(\phi(0)=0)$ and {\bf (P2)}$\Leftrightarrow (\phi(1)=0)$. 
In the case of the $C_\theta$ function, $\Delta$ rewrites
$$
\Delta(u_1, u_2, v_1, v_2)=(u_2-u_1)(v_2-v_1) + \theta(\phi(u_2)-\phi(u_1))(\phi(v_2)-\phi(v_1)),
$$
and then {\bf (b)}$\Rightarrow$ {\bf (P3)}.
Conversely, if there are  $x<y$ such that $|\phi(x)-\phi(y)|>|x-y|$,
then with the choice $\theta=-1$, $u_1=v_1=x$ and $u_2=v_2=y$,  
we obtain:\\ $ \Delta(x,y,x,y)=(x-y)^2 - (\phi(x)-\phi(y))^2 < 0.$
Therefore {\bf (P3)} $\Rightarrow$ {\bf (b)}.
Lastly, $C_\theta$ is absolutely continuous since $\phi$ satisfies the Lipschitz condition.
\hfill $\Box$\\

\noindent
The following corollary provides a new characterization of the functions $\phi$ 
generating parametric families of copulas.
It will be of great help in the following.  
\begin{coro}
\label{corocop}
$\phi$ generates a parametric family of copulas $\{C_\theta,\; \theta\in[-1,1]\}$
 if and only if it satisfies the following conditions~:
\begin{description}
\item [(i)] $\phi$ est absolutely continuous,
\item [(ii)] $|\phi'(x)|\leq 1$ almost everywhere in $I$,
\item [(iii)] $|\phi(x)|\leq \min(x,1-x)$, $\forall x\in I$.
\end{description}
In such a case, $C_\theta$ is absolutely continuous.
\end{coro}
\noindent {\bf Proof: } From a classical analysis result (\cite{ANALYSE}, Lemma 2.1), {\bf (b)} is satisfied if and only if {\bf (i)} and {\bf (ii)} both hold. Now, assume that conditions {\bf (i)--(iii)} are satisfied. Taking successively $x=0$ and $x=1$ in {\bf (iii)} yields {\bf (a)}. Therefore the corollary implies the theorem.\\
Conversely, assume that {\bf (b)} holds. Taking successively $y=0$ and $y=1$ in the Lipschitz condition gives {\bf (iii)}. So, the theorem implies the corollary.
\hfill $\Box$\\

\noindent
The function $\phi$ plays a role similar to the generating function 
in Archimedian copulas \cite{GENEST}. Each copula $C_\theta$ is entirely described by the univariate function $\phi$ and the parameter $\theta$, which tunes the dependence between the margins (see Section \ref{secdep}). Symmetry and dependence properties of the copula $C_\theta$ will found a geometrical interpretation on the graph of $\phi$. Furthermore, the choice of $\phi$ determines the vertical and horizontal sections of the copula up to a multiplicative factor and an additive linear function.\\
\noindent  We now give some examples of functions $\phi$ generating parametric families of copulas. These examples will be used to illustrate the symmetry and dependence properties introduced in the following sections.

\begin{ex} The following functions $\phi$ generate parametric families of copulas.
\label{ex1}
\label{exemples}
\begin{enumerate}
\item $\phi^{[1]}(x)=\min(x,1-x)$ is the upper bound of Corollary 
\ref{corocop}{\bf (iii)}. 
\item $\phi^{[2]}(x)=x(1-x)$ generates the Farlie-Gumbel-Morgenstern (FGM) family of copulas \cite{FGM}, which contains all copulas with both horizontal and vertical quadratic sections \cite{QUESADA}. 
\item $\phi^{[3]}(x)=x(1-x)(1-2x)$ defines the parametric family of symmetric copulas with cubic sections proposed in \cite{NELSEN97}, equation (4.4). 
\item $\phi^{[4]}(x)=\frac{1}{\pi}\sin(\pi x)$ induces a family of copulas able to modelize large dependences. 
\end{enumerate}
\end{ex}
In the sequel, we note $\{C_\theta^{[i]}\}$ the parametric family of copulas
associated to the function $\phi^{[i]}$, $i\in\{1,\dots,4\}$.
Graphes of the functions $\phi^{[i]}$ are plotted in~\cite{Nous}.
As a consequence of the condition {\bf (iii)} of Corollary \ref{corocop}, graphes of $\phi^{[1]}$ and $(-\phi^{[1]})$ are the edges of a square $K$, within which lies the graph of $\phi$.

\noindent In the following, we review several concepts of symmetry and dependence.
Throughout the rest of this paper, $\phi$ denotes a function satisfying 
the conditions given in Theorem \ref{theocop}, and $\{C_\theta\}$ represents the parametric family of copulas genetated by the function $\phi$.

\section{Symmetry properties}
\label{secsym}

Let $(a,b)\in\mathbb{R}^2$ and $(X,Y)$ a random pair. We say that $X$ is symmetric about $a$ if the cumulative distribution function of $(X-a)$ and $(a-X)$ are identical. The following definitions generalize this symmetry concept to the bivariate case: 
\begin{itemize} 
\item $X$ and $Y$ are exchangeable if $(X,Y)$ and $(Y,X)$ are identically
distributed.
\item $(X,Y)$ is marginally symmetric about $(a,b)$ if $X$ and $Y$ are symmetric about $a$ and $b$ respectively.
\item $(X,Y)$ is radially symmetric about $(a,b)$ if $(X-a,Y-b)$ and $(a-X,b-Y)$ follow the same joint cumulative distribution function.
\item $(X,Y)$ is jointly symmetric about $(a,b)$ if the pairs $(X-a,Y-b)$, $(a-X,b-Y)$, $(X-a,b-Y)$ and $(a-X,Y-b)$ have a common joint cumulative distribution function.
\end{itemize}
The following theorem provides conditions on $\phi$ to ensure that the couple $(X,Y)$ with  associated copula $C_\theta$ is radially (or jointly) symmetric.
\begin{theo}
\label{thsym}
\hfill
\begin{description}
\item [(i)] If $X$ and $Y$ are identically distributed then $X$ and $Y$ are exchangeable.
\end{description}
Besides, if $(X,Y)$ is marginally symmetric about $(a,b)$ then:
\begin{description}
\item [(ii)] $(X,Y)$ is radially symmetric about $(a,b)$ if and only if 
\begin{equation}
\label{eqsym}
\mbox{ either } \forall u \in I,\; \phi(u)=\phi(1-u) \mbox{ or }\forall u \in I,\; \phi(u)=-\phi(1-u).
\end{equation}
\item [(iii)] $(X,Y)$ is jointly symmetric about $(a,b)$ if and only if $\forall u \in I,\; \phi(u)=-\phi(1-u)$.
\end{description}
\end{theo}
\noindent {\bf Proof: } 
\begin{description}
\item [(i)] When $X$ and $Y$ are identically distributed, exchangeability
is equivalent to the symmetry of the copula. In other words, $C_\theta$ must 
verify $C_\theta(u,v)=C_\theta(v,u)$,
$\forall (u,v)\in I^2$, which is the case by definition, see (\ref{defcop}).
\item [(ii)]  Assume $(X,Y)$ is marginally symmetric about $(a,b)$.
Then, from  Theorem 2.7.3 in \cite{NELSEN99}, $(X,Y)$ is radially symmetric about $(a,b)$ if and only if 
\begin{equation}
\label{eqdelta}
\forall (u,v)\in I^2,\;\; \delta_\theta(u,v)=0,
\end{equation}
where we have defined
$$
\delta_\theta(u,v)=C_\theta(u,v)-C_\theta(1-u,1-v)-u-v+1
=\theta[\phi(u)\phi(v)-\phi(1-u)\phi(1-v)].
$$
It clearly appears that (\ref{eqsym}) implies (\ref{eqdelta}).
Conversely, suppose that (\ref{eqdelta}) is verified.
Then, $\delta_1(u,u)=\phi^2(u)-\phi^2(1-u)=0$ for all $u\in I$ and
consequently, either $\phi(u)=\phi(1-u)$ or $\phi(u)=-\phi(1-u)$. 
If there exist $(u_1,u_2)$ such that $\phi(u_1)=\phi(1-u_1)\neq 0$ and $\phi(u_2)=-\phi(1-u_2)\neq 0$, then $\delta_1(u_1,u_2)=2\phi(u_1)\phi(u_2)\neq 0$. 
As a conclusion (\ref{eqdelta}) implies (\ref{eqsym}).
\item [(iii)] Assume $(X,Y)$ is marginally symmetric about $(a,b)$.
$(X,Y)$ is jointly symmetric about $(a,b)$ if and only if
\begin{equation}
\label{eqdelta2}
\forall(u,v) \in I^2,\;\; \delta'_\theta(u,v)=0,
\end{equation}
where we have defined
$$
\delta'_\theta(u,v)=C_\theta(u,v)+C_\theta(u,1-v)-u=\theta\phi(u)[\phi(v)+\phi(1-v)].
$$
It appears immediately that (\ref{eqdelta2}) is equivalent to
$\phi(u)=-\phi(1-u)$, $\forall u \in I$.
\hfill $\Box$
\end{description}
As an example, any marginally symmetric random pair $(X,Y)$ associated to 
a copula of the families $\{C_\theta^{[i]}\}$, $i\in\{1,\dots,4\}$ is radially
 symmetric. Moreover, in the case $i=3$, $(X,Y)$ is jointly symmetric.
We now focus on the dependence properties of the semiparametric family of copulas.

          \section{Measures of association}
\label{secdep}
In the next two sections, we note $(X,Y)$ a random pair with joint distribution $H$, density distribution $h$, copula $C$ and margins $F$ and $G$. The case $C=C_\theta$ will be explicitly precised. 
Three invariant to strictly increasing function measures of association between the components of the random pair $(X,Y)$
are usually considered:
\begin{itemize}
\item the normalized volume between graphes of $H$ and $FG$ \cite{SIGMA}, $$\sigma=12\int_0^1\!\!\! \int_0^1 |C(u,v)\!-\!uv|dudv,$$
\item the Kendall's Tau [2,3], defined as the probability of concordance minus the probability of discordance of two pairs $(X_1,Y_1)$ and $(X_2,Y_2)$ described by the same joint bivariate law $H$,
\begin{equation}
\label{coefftau}
\tau=4\int_0^1 \!\!\! \int_0^1 C(u,v) dC(u,v)\!-\!1,\;
\end{equation}
\item the Spearman's Rho [2,3], which is the probability of concordance minus the probability of discordance of two pairs $(X_1,Y_1)$ and $(X_2,Y_2)$ with respective joint cumulative law $H$ and $FG$, $$\rho=12\int_0^1 \!\!\! \int_0^1 C(u,v) dudv\!-\!3.$$
\end{itemize}
In the case of a copula generated by (\ref{defcop}),
these measures rewrite only in terms of the function~$\phi$.
\begin{prop}
Let $(X,Y)$ be a random pair with copula $C_\theta$ given by (\ref{defcop}).
The association coefficients are:
$$
\sigma_\theta=12|\theta| \left(\int_0^1 |\phi(u)|du\right)^2,\;\;
\tau_\theta=8\theta \left(\int_0^1 \phi(u)du\right)^2,\;\;
\rho_\theta=12\theta \left(\int_0^1 \phi(u)du\right)^2=\frac{3}{2}\tau_\theta.
$$
\end{prop}
\noindent {\bf Proof:} The proof are very similar for the three coefficients.
Let us consider the example of the Kendall's Tau.
The density distribution of the cumulative distribution function $C_\theta$ is 
\begin{equation}
\label{eqdens}
c_\theta(u,v)=1+\theta \phi'(u) \phi'(v).
\end{equation}
Replacing in (\ref{coefftau}), it yields 
\begin{eqnarray*}
\tau_\theta & = & 4\theta \left[\int_0^1 \phi(u)du\right]^2+
		  4\theta \left [\int_0^1 u\phi'(u)du\right]^2+
		  4\theta^2 \left[\int_0^1 \phi(u)\phi'(u)du\right]^2,\\
            & = & 4\theta \left[\int_0^1 \phi(u)du\right]^2+
		  4\theta \left[\phi(1)-\int_0^1 \phi(u)du\right]^2+
		  2\theta^2[\phi^2(1)-\phi^2(0)],
\end{eqnarray*}
after partial integration. Theorem \ref{theocop} entails $\phi(1)=\phi(0)=0$ and the conclusion
follows.
\hfill $\Box$\\

\noindent The measures $\tau_\theta$ and $\rho_\theta$ linearly increase with $\theta$, which appears as an association parameter. They are also proportionnals to the square of the surface lying between the graph of $\phi$ and the $x$-axis.  
Bounds for each of these measures of association are deduced from
 Corollary~1~{\bf (iii)}.
\begin{prop}
\label{bounds}
$\forall \theta\in [-1,1]$, $0\leq\sigma_\theta\leq 3|\theta|/4$,
$|\tau_\theta|\leq |\theta|/2$ and $|\rho_\theta|\leq 3|\theta|/4$.
\end{prop}
Therefore the range of each association coefficient is
$0\leq\sigma_\theta\leq 3/4$, $-1/2\leq\tau_\theta\leq 1/2$
and $-3/4\leq\rho_\theta\leq 3/4$.
The semiparametric family of copulas defined by (\ref{defcop}) seems to be a good tool to build low complexity models with moderate dependences ($0\leq |\tau| \leq 1/2$). This is illustrated in Table \ref{tabcoeff} where the values obtained for the copulas presented in Example 1
are reported. 
Results obtained for the cubic sections copulas illustrate the fact that $\tau$ and $\rho$ imperfectly measure the association relationships.
 All the values associated to $\{C^{[4]}_\theta\}$ are larger than those obtained with the FGM family. 
Lastly, bounds of Proposition \ref{bounds} are reached with family $\{C^{[1]}_\theta\}$. This hierarchy of results graphically appears in~\cite{Nous}.
For each fixed $\theta$, the closer to the edges of the square $K$ is the graph of $\phi$, the larger is the association coefficient of the copula induced by $\phi$.
Non-derivability of $\phi^{[1]}$ can be a drawback. In order to build  more regular copulas, with measures of association remaining closed from optimal values, we define the
sequence of ${\cal C}^1$ functions $(\phi_n^{[5]})$ as following: 
$$
\forall n>0,\; \phi_n^{[5]}(x)=\left|
\begin{array}{lll}
\vspace*{1mm}
x & \mbox{ if} & 0\leq x \leq \displaystyle \frac{1}{2}-\frac{1}{n} \\
\vspace*{1mm}
\displaystyle \frac{1}{2}\left(1-\frac{1}{n}\right)-\frac{n}{2}\left(x-\frac{1}{2}\right)^2 & \mbox{ if } 
		& \displaystyle \frac{1}{2}-\frac{1}{n}< x < \frac{1}{2}+\frac{1}{n} \\
\vspace*{1mm}
1-x & \mbox{ if } & \displaystyle \frac{1}{2}+\frac{1}{n}\leq x \leq 1.
\end{array}
\right.
$$
Moreover each copula is derivable and the associated Kendall's Tau sequence is given
 by 
$$
\tau^{[5]}_{n,\theta}=8\theta\left(\frac{1}{4}-\frac{1}{3n^2}\right)^2, \;\; n>0.
$$ 
It converges to the optimal Tau, $\tau^{[5]}_{n,\theta}\to \theta/2$ as $n\to\infty$. 
It is also possible to build ${\cal C}^{\infty}$ copulas, with measures of association closed from optimal values. Define the sequence $(\phi_n^{[6]})$ of ${\cal C}^{\infty}$ functions by:
$$
\forall n\geq 2, \;\;\phi_n^{[6]}(x)=1-\left(x^n+(1-x)^n\right)^{1/n}.
$$
Similarly to the previous example, 
for each $n\geq 2$, $\phi_n^{[6]}$ induces a parametric family of copulas $\{C^{[6]}_{n,\theta}\}$.
Moreover $\forall x\in I,\; \phi_n^{[6]}(x)\to\phi^{[1]}(x)$ as $n\to\infty$ and $|\phi_n^{[6]}|\leq 1$.
So, by Lebesgue dominated convergence theorem, the associated Kendall's Tau, $\tau^{[6]}_{n,\theta}$ converges to the optimal value as $n\to\infty$.
Graphes of $\phi_2^{[6]}$, $\phi_4^{[6]}$, and $\phi_8^{[6]}$ are plotted in~\cite{Nous}.
 We clearly see that the surface lying between the $\phi_n^{[6]}$ graph and 
the $x$-axis increases with $n$. The sequences $(\tau^{[6]}_{n,1})$ and 
$(\tau^{[5]}_{n,1})$ are plotted in~\cite{Nous}. We can observe that for $n\geq 2,\;\; \tau^{[5]}_{n,1}$ and $\tau^{[6]}_{n,1}$ are both larger than the value obtained by the FGM copula $C^{[2]}_1$ and that they exceed the value reached by the copula $C^{[4]}_1$, whenever $n\geq 3$.

\section{Concepts of dependence}
\label{secdep2}

Although the notion of independence between  $X$ and $Y$ is clearly defined by $H=FG$, the concept of dependence may be defined in different ways. A dependence ordering,
indicating if the random pair $(X_1,Y_1)$ of margins $F$ and $G$ is ``more dependent'' than the pair $(X_2,Y_2)$ with the same margins,
 can be associated to each definition.
We firstly study the properties of positive dependence of a copula $C_\theta$ and secondly we search to order elements within the parametric family $\{C_\theta\}$.

\subsection{Positive dependence}

In this subsection, we will use several concepts of positive dependence, which express that two variables are large (or small) simultaneously. 
Let us shortly review their definition.

\begin{itemize}
\item PFD: $X$ and $Y$ are Positive Function Dependent if for all integrable real-valued function $g$ $\mathbb{E}_h[g(X)g(Y)]-\mathbb{E}_h[g(X)]\mathbb{E}_h[g(Y)]\geq 0$ where $\mathbb{E}_h$ is the expectation symbol relative to the density $h$.
\item PQD: $X$ and $Y$ are Positively Quadrant Dependent if 
\begin{equation}
\label{PQD}
\forall (x,y)\in \mathbb{R}^2,\;\; P(X\leq x,Y\leq y)\geq P(X\leq x)P(Y\leq y).
\end{equation}
\item  LTD$(Y|X)$: $Y$ is Left Tail Decreasing in $X$ if $P(Y\leq y| X\leq x)$ is nonincreasing in $x$ for all $y$. A similar definition can be given for $LTD(X|Y)$.
\item RTI$(Y|X)$: $Y$ is Right Tail Increasing in $X$ if $P(Y>y|X>x)$ is nondecreasing in $x$ for all~$y$.
\item SI$(Y|X)$: $Y$ is Stochastically Increasing in $X$ if $P(Y>y|X=x)$ is nondecreasing in $x$ for all~$y$.
\item LCSD: $X$ and $Y$ are Left Corner Set Decreasing  if $P(X \leq x,Y \leq y|X \leq x',Y \leq y')$ is nonincreasing in $x'$ and $y'$ for all $x$ and $y$.
\item RCSI: $X$ and $Y$ are Right Corner Set Increasing  if $P(X>x,Y >y|X>x',Y>y')$ is nondecreasing in $x'$ and $y'$ for all $x$ and $y$.
\item TP2 density: $(X,Y)$ have the TP2 density property if $h$ is a totally positive function of order 2 {\it i.e.}
$h(x_1,y_1)h(x_2,y_2)-h(x_1,y_2)h(x_2,y_1)\geq 0$ for all $(x_1,x_2,y_1,y_2)\in 
I^4$ such that $x_1\leq x_2$ and $y_1\leq y_2$.
\end{itemize}
When $X$ and $Y$ are exchangeable, there are no reason to distinguish SI$(Y|X)$ and SI$(X|Y)$, which will be both noted SI. Similarly, we will denote LTD the equivalent properties LTD$(Y|X)$ and LTD$(X|Y)$, and RTI, RTI$(Y|X)$ or RTI$(X|Y)$.
The links between these concepts \cite{LEHANN,NELSEN91} are illustrated
in~\cite{Nous}.
The following theorem is devoted to the study of properties of positive dependence of any pair $(X,Y)$ associated with the copula $C_\theta$ defined by  (\ref{defcop}). Similar results can be established for the corresponding concepts of negative dependence. 
\begin{theo}
\label{theodep1}
Let $\theta>0$ and $(X,Y)$ a random pair with copula $C_\theta$.
\begin{description}
\item [(i)] $X$ and $Y$ are PFD.
\item [(ii)] $X$ and $Y$ are PQD if and only if either $\forall u \in I,\; \phi(u)\geq 0$  or $\forall u \in I,\; \phi(u)\leq 0$.
\item [(iii)] $X$ and $Y$ are LTD if and only if $\phi(u)/u$ is monotone.
\item [(iv)] $X$ and $Y$ are RTI if and only if $\phi(u)/(u-1)$ is monotone.
\item [(v)] $X$ and $Y$ are LCSD if and only if they are LTD.
\item [(vi)] $X$ and $Y$ are RCSI if and only if they are RTI.
\item [(vii)] $X$ and $Y$ are SI  if and only if $\phi$ is either concave or convex.
\item [(viii)] $X$ and $Y$ have the TP2 density property if and only if 
they are SI.
\end{description}
\end{theo}
\noindent {\bf Proof: }
\begin{description}
\item [(i)] Let $g$ be an integrable real-valued function on $I$.
The density distribution $c_\theta$ of the cumulative distribution $C_\theta$ 
is given by (\ref{eqdens}). Routine calculations yield
$$
\mathbb{E}_{c_\theta}(g(X)g(Y))-\mathbb{E}_{c_\theta}(g(X))
\mathbb{E}_{c_\theta}(g(Y))=\theta {\left[\int_0^1g(t)\phi'(t)dt\right]}^2 \geq 0,
$$
since $\theta\geq 0$.
\item [(ii)] The pair $(X,Y)$ is PQD if and only if the uniform I-margins pair $(U,V)$ with distribution $C_\theta$ is PQD. For $(U,V)$, condition (\ref{PQD}) simply rewrites $\theta\phi(u)\phi(v)\geq 0$,
$\forall (u,v)\in I^2$ and the conclusion follows.
\item[(iii)] Since $C_{\theta}$ is symmetric, the necessary and sufficient
conditions given by Theorem 5.2.5 in \cite{NELSEN99}, simply reduce to the unique condition: 
$C_\theta(u,v)/u=v+\theta\phi(v)\phi(u)/u$
is nonincreasing in $u$ for all $v\in I$.
Suppose for instance $\phi(u)/u$ is nonincreasing. Then, since $\phi(1)=0$,
we have $\phi(v)\geq 0$ for all $v$ in $I$ and thus $C_\theta(u,v)/u$
is nonincreasing in $u$ for all $v\in I$. The case of a nondecreasing function $\phi(u)/u$ is similar.\\
Conversely, suppose $C_\theta(u,v)/u$ is nonincreasing in $u$ for all $v\in I$.
If there exist $v\in I$ such that $\phi(v)> 0$ then $\phi(u)/u$ is
nonincreasing too. The case where there is $v\in I$ such that $\phi(v)<0$ is similar.
The case $\phi(v)=0$ for all $v\in I$ is trivial.
As a conclusion, we have shown that $C_\theta(u,v)/u$ is nonincreasing in $u$ for all $v\in I$ if and only if $\phi(u)/u$ is monotone.
\item [(iv)] is similar to {\bf (iii)}.
\item [(v)] In view of  Corollary 5.2.17 in \cite{NELSEN99}, a necessary and sufficient
condition for $X$ and $Y$ LCSD is: $C_\theta$ is a totally positive function
of order 2.
For all $u_1 \leq u_2$ and $v_1\leq v_2$, the quantity 
$$
C_\theta(u_1,v_1)C_\theta(u_2,v_2)-C_\theta(u_1,v_2)C_\theta(u_2,v_1)=
\theta\left[\frac{\phi(v_2)}{v_2}-\frac{\phi(v_1)}{v_1}\right]\left[\frac{\phi(u_2)}{u_2}-\frac{\phi(u_1)}{u_1}\right]u_1 u_2 v_1 v_2
$$
is nonnegative if and only if $\phi(u)/u$ is monotone.
This is the necessary and sufficient condition for $X$ and $Y$ LTD given in {\bf (iii)}.
\item [(vi)] is similar to {\bf (v)}.
\item [(vii)] In view of the $C_\theta$ symmetry, the geometric interpretation
of stochastic monotonicity given by Corollary 5.2.11 in \cite{NELSEN99} provides
the necessary and sufficient condition:
$C_\theta(u,v)$ is a concave function of $u$.
Suppose $\phi$ is a concave function. Then, taking into account that
$\phi(0)=\phi(1)=1$, we easily show that $\phi(v)\geq 0$ for all $v\in I$
and therefore $C_\theta(u,v)$ is a concave function of $u$.
The case of a convex function $\phi$ is similar.\\
Conversely, suppose $C_\theta(u,v)$ is a concave function of $u$.
If there exist $v$ such that $\phi(v)>0$ then (\ref{defcop}) shows
that $\phi$ is concave. The case where there exist $v\in I$ such that
$\phi(v)<0$ is similar.
The case $\phi(v)=0$ for all $v\in I$ is trivial.
As a conclusion we have shown the equivalence between 
$C_\theta(u,v)$ is a concave function of $u$ and $\phi$ is concave or convex.
\item[(viii)]
$X$ and $Y$ have the TP2 density property if and only if
the density of the copula verifies
$$
\forall u_1 \leq u_2,\;v_1\leq v_2,\;\; c_\theta(u_1,v_1)c_\theta(u_2,v_2)-c_\theta(u_1,v_2)c_\theta(u_2,v_1) \geq 0,
$$
which rewrites
$$
\forall u_1 \leq u_2,\;v_1\leq v_2,\;\; [\phi'(u_1)-\phi'(u_2)][\phi'(v_1)-\phi'(v_2)]\geq 0.
$$
Equivalently, $\phi'$ is either nonincreasing or nondecreasing
which means that $\phi$ is either concave or convex.
This is the necessary and sufficient condition for $X$ and $Y$ SI established
in {\bf (vii)}.
\hfill $\Box$
\end{description}
Results of Theorem \ref{theodep1} are illustrated in~\cite{Nous}.
For instance, random pairs with copula
in the parametric families generated by $\phi^{[1]}$, $\phi^{[2]}$,
$\phi^{[4]}$, $\phi_n^{[5]}$ or $\phi_n^{[6]}$ have the TP2 density property
since these functions are concave. 
Consequently, these random pairs are also SI, LCSD, RCSI, LTD, RTI, PQD, and PFD.  When the copula belongs to the family $\{C_\theta^{[3]}\}$, 
the associated random pair is PFD but not PQD. \\
Finally, let us note that
Theorem~\ref{theodep1} provides a very simple tool to prove that
the family proposed in \cite{LAI}:
$
c(u,v)=uv+\theta u^b v^b(1-u)^a (1-v)^a,\; \theta\in[-1,1],\;a\geq 1,\;b\geq 1
$
is PQD.

\subsection{Dependence orderings}

In Theorem \ref{theodep1} positive dependence properties of the
copula $C_\theta$ are established for fixed $\theta>0$.
In this subsection we investigate the variations of
the ``dependence amount'' when $\theta$ increases.
To this end, it is necessary to define a dependence ordering
which can decide whether
a multivariate cumulative distribution function is more
dependent than another.
Consider $H$ and ${H^*}$ two bivariate cumulative
distribution functions with the same margins.
We study two dependence orderings:
\begin{itemize}
\item ${H^*}$ is more PQD (or more concordant) than $H$ if
$H(x_1,x_2)\leq {H^*}(x_1,x_2)$ for all $(x_1,x_2)$.
We note $H \prec_C {H^*}$.
\item Note $H_{2|1}$ and ${H^*}_{2|1}$ the conditional
 distributions of the second random variable given the first one.
${H^*}$ is more SI than $H$ if
${H^*}_{2|1}^{-1}(H_{2|1}(x_2|x_1)|x_1)$ is increasing in $x_1$.
We note $H \prec_{SI} {H^*}$.
\end{itemize}
Remark that $H \prec_{SI} {H^*}$ 
implies $H \prec_C {H^*}$ (see Theorem 2.12 in \cite{JOE}).
The natural extension of these definitions to parametric families of copulas
is the following.
\begin{itemize}
\item A family of copulas $\{C_\theta\}$ is ordered in concordance if
$\theta \leq \theta'$ implies $C_\theta \prec_C C_{\theta'}$.
\item A family of copulas $\{C_\theta\}$ is SI ordered if
$\theta \leq \theta'$ implies $C_\theta \prec_{SI} C_{\theta'}$.
\end{itemize} 
As mentionned in \cite{JOE}, $C_\theta \prec_C C_{\theta'}$ implies $\tau_{\theta}
\leq \tau_{\theta'}$ and $\rho_{\theta}\leq \rho_{\theta'}$.
We now give conditions under which
the parametric family of copulas obtained by choosing $\phi$
according to Corollary \ref{corocop} is
stochasticly ordered.
\begin{theo}\hfill
\begin{description}
\item [(i)] $\{C_\theta\}$ is ordered in concordance if and only if $\phi$ is either positive or negative.
\item [(ii)] Assume $\phi$ is 2-derivable.
$\{C_\theta\}$ is SI ordered if $\phi$ is either convex or concave.
\end{description}
\end{theo}
\noindent {\bf Proof: }
\begin{description}
\item [(i)] The proof is straigtforward.
\item [(ii)] A direct proof is not possible since, in the general case,
the inverse of the conditional cumulative distribution function does not have a  
closed-form. A necessary condition is given by Theorem 2.14 in \cite{JOE}:
$\{C_\theta\}$ is SI ordered if for all increasing real-valued
function $g$,
$$
A(u,v,\theta)=\frac{\partial^2 B}{\partial v \partial u}\frac{\partial B}{\partial \theta}-\frac{\partial^2 B}{\partial \theta \partial u}\frac{\partial B}{\partial v} \geq 0,
$$
where $B(u,v,\theta)= g(v+\theta \phi'(u)\phi(v))$.
We obtain 
$$
A(u,v,\theta)= -\phi''(u)\phi(v){[g'(v+\theta\phi'(u)\phi(v))]}^2.
$$
The sign of $A(u,v,\theta)$ is given by $-\phi''(u)\phi(v)$, which is positive 
whenever $\phi$ is concave or convex.
\hfill $\Box$\\
\end{description}
For instance, all families discussed in this document are ordered in concordance,
but $\{C_\theta^{[3]}\}$ is not SI ordered.

\section{Conclusion and further work}
In this paper, a symmetric semiparametric family of copulas is studied
and its symmetry and dependence properties are established.
As for Archimedian copulas, a parametric subfamily is generated by an
 univariate function and both symmetry and dependence properties of the family can be geometrically interpreted on the graph of the generating function. 
Furthermore, the horizontal and vertical sections of the copula
are essentially described by the generating function.
Finally, the proposed family of copulas provides a ``simple'' way to modelize moderate degrees of dependence, in the sense that the copulas  are convenient to perform classical calculations on a probability law. 
Our further work will consist in generalizing the family definition to the asymmetric case:
\begin{equation} 
\label{asym}
C_\theta(u,v)=uv+\theta a(u)b(v),\; \theta\in[-1,1],
\end{equation} 
where $a$ and $b$ are functions defined on $I$. This family was introduced
in \cite{FLO}, chapter 4. It can be shown that, 
similarly to Theorem \ref{theocop},
necessary and sufficient conditions for (\ref{asym}) to be a copula are
$a(0)=a(1)=b(0)=b(1)=0$ and 
$|a(x)-a(y)| |b(z)-b(t)| \leq |x-y| |z-t|$ for all $(x,y,z,t)\in I^4$.
The study of the related measures of association and dependence concepts
will be the topic of a next paper.

\section*{Acknowledgements}

This work has been done while the second author was visiting the
Centre de Recherche Math\'e\-ma\-ti\-ques, Universit\'e de Montr\'eal,
whose hospitality is gratefully acknowledged.


\begin{table}[h]
\begin{center}
\begin{tabular}{|c|c|c|c|}
\hline
$\phi(x)$ & $\sigma_\theta$ & $\tau_\theta$ & $\rho_\theta$ \\
\hline
&&&\\
$\phi^{[1]}(x)=\min(x,1-x)$ & $3|\theta|/4$ & $\theta/2$  & $3\theta/4$ \\
$\phi^{[2]}(x)=x(1-x)$      & $|\theta|/3$  & $2\theta/9$ & $\theta/3$ \\
$\phi^{[3]}(x)=x(1-x)(1-2x)$& $3|\theta|/64$& 0           & 0 \\
$\phi^{[4]}(x)=\frac{1}{\pi}\sin{(\pi x)}$ & $48|\theta|/{\pi^4}$ & $32\theta/{\pi^4}$ &
 $48\theta/{\pi^4}$ \\
&&&\\
\hline
\end{tabular}
\end{center}
\caption{Measures of dependence associated to the Example 1 copulas.}
\label{tabcoeff}
\end{table}
\end{document}